\documentclass{amsart}

\usepackage{graphicx}
\usepackage{amssymb}

\newtheorem{defn}{Definition}

\newtheorem{rem}{Remark}

\begin{document}

\title[Chebyshev Series of the Feigenbaum Function]{Chebyshev Series Representation of Feigenbaum's Period-Doubling Function}

\author{Richard J. Mathar}
\urladdr{http://www.strw.leidenuniv.nl/~mathar}
\email{mathar@strw.leidenuniv.nl}
\address{Leiden Observatory, Leiden University, P.O. Box 9513, 2300 RA Leiden, The Netherlands}

\subjclass[2010]{Primary 37G15, 37-04; Secondary 37M20}

\date{\today}
\keywords{Feigenbaum constant, period doubling, Chebyshev series}

\begin{abstract}
The Feigenbaum-Cvitanovi\'{c}
equation $-\lambda g(x)=g(g(\lambda x))$ is solved over
the interval $0\le x\le 1$
with a Chebyshev series representation of $g(x)$.
Accurate expansion coefficients are tabulated for solutions 
$g(x)=1+O(x^z)$
with even exponents from $z=2$ up to $z=14$.
\end{abstract}

\maketitle
\section{Introduction}

This article provides high precision approximations of functions $g(x)$
which solve the Feigenbaum--Cvitanovi\'c equation
\cite{FeigenbaumJSP19,FeigenbaumJSP21,ColletCMP76,TsygCRAS334,EpsteinCMP106,EpsteinCMP81,ChristiaJPA23}
\begin{equation}
-\lambda  g(x)=g(g(\lambda x)),
\label{eq.f}
\end{equation}
scaled such that
\begin{equation}
g(0)=1.
\label{eq.gspec}
\end{equation}
The parameter $\lambda$
plays the role of an eigenvalue
bound to the solutions via
\begin{equation}
g(1) = -\lambda.
\end{equation}
Further below we shall refer to $1/\lambda$ as \emph{the} Feigenbaum
constant(s)---although in
a broader context
a larger variety of numbers carries that name.

Only even functions $g(x)=g(-x)$ are discussed, so the standard
representation is the Taylor series
\cite{FeigenbaumJSP21,BriggsPhd,BriggsMC57}.
\begin{defn} (Taylor series coefficients $b_n$)
\begin{equation}
g(x) = 1+\sum_{n=z,2z,3z,\ldots} b_n x^n;\quad z=2,4,6,\ldots
\label{eq.gtay}
\end{equation}
\end{defn}

In this manuscript, the function $g(x)$ is expanded in a series of Chebyshev
Polynomials $T(x)$ \cite[\S 22]{AS}\cite[\S 18]{DLMF}, which---for well understood
reasons---supplies a more stable basis
than the bare powers \cite{Cody}.
\begin{defn} (Chebyshev series expansion coefficients $t_n$).
\begin{equation}
g(x)=\sideset{}{'}\sum_{n\ge 0} t_n T_n(x^d);\quad d=1,2,3,\ldots .
\label{eq.gx}
\end{equation}
The prime at the sum symbol indicates
the term at $n=0$ is halved.
\end{defn}
We are considering even exponents $z$,
so $g(x)$ is even and $t_n=0$ whenever $n$ is odd.
The integer variable $d$
plays the following role: Linear coupling equations
between the $t_n$ would have to be set up to remove the contributions by
powers $x^n$ ($z\nmid n$) from the solutions
if $z>2$ and $d=1$.
(The origin of this constraint that all exponents of $x$ be multiples of 
the same $z$ is not elaborated here.)
Explicitly, equating equal powers of $x$ in (\ref{eq.gtay}) and (\ref{eq.gx}),
the non-zero values of $b_n$ are
\begin{equation}
2^{n/d}\sideset{}{'}
\sum_{\substack{s\ge n/d\\ s+n/d \equiv 0 \pmod 2}}
t_s
\frac{s}{s+n/d}(-1)^{(s-n/d)/2}\binom{(s+n/d)/2}{n/d}
= b_n
,\quad d\mid n
.
\label{eq.boft}
\end{equation}
In the algorithm described further down, these would
multiply the order of the linear system of equations by an
approximate factor of $z/2$.
For enhanced efficiency, $d$ will be taken as
\begin{equation}
d=z/2,
\end{equation}
which ensures that only powers of $x$ with exponents
divisible by $z$ enter the calculation.

Equation (\ref{eq.gspec}) and the special values of the Chebyshev polynomials
induce a sum  rule
and a
coupling with $\lambda$:
\begin{gather}
g(0)=1\therefore 
1=
\sideset{}{'}\sum_n
(-1)^{n/2}
t_n
;
\label{eq.g0}
\\
g(1)= -\lambda =
\sideset{}{'}\sum_n t_n
.
\label{eq.g1}
\end{gather}

\section{Numerical Algorithm}
\subsection{Multivariate Newton Iteration}

The right hand side of (\ref{eq.f}) is
\begin{gather}
g(g(\lambda x))
= 
\sideset{}{'}\sum_m t_m T_m\left(g^d[\lambda x]\right)
=
\sideset{}{'}\sum_m t_m T_m\left(\left\{\sideset{}{'}\sum_n t_n T_n(\lambda^d x^d)\right\}^d\right)
.
\end{gather}
Moving all terms to one side of the equation,
the aim is to obtain a zero of the function $f$,
\begin{gather}
f_j(\{t_m\})\equiv 
\sideset{}{'}\sum_m t_m T_m\left(\left\{\sideset{}{'}\sum_n t_n T_n(\lambda^d x_j^d)\right\}^d\right)
+\lambda
\sideset{}{'}\sum_n t_n T_n(x_j^d) =0;
\quad j=1,2,\ldots
\label{eq.fzer}
\end{gather}
We solve for $g$ with a finite, iteratively enlarged set of expansion
coefficients 
$t_n$ which are obtained by fitting at a finite set of 
the standard Chebyshev abscissa
points
$x_j=\cos\theta_j$, $\theta_j= j\pi/N$.
\begin{rem}
With the presence of $d$, which effectively replaces the Chebyshev weights
$1/\sqrt{1-x^2}$ by $x^{d-1}/\sqrt{1-x^{2d}}$, these abscissae may not be the
optimum choice if $d>1$.
\end{rem}
The common multivariate Newton algorithm with
a $N\times N$ matrix of first derivatives is employed:
we start with a set of approximations $t_k$,
and calculate a vector of corrections $\Delta t_k$
which are the solutions to the linear system of equations
\begin{gather}
f_j(\{t_m\})
+\sum_k \frac{\partial f_j}{\partial t_k}\Delta t_k =0;
\quad j=1,2,\ldots 
\end{gather}
This is actually done on $N-1$ abscissa points $x_j$, because 
one row of the system of equations is reserved to accommodate
(\ref{eq.g0}):
\begin{equation}
\frac12 t_0-t_2+t_4-\ldots -1=0.
\end{equation}
\begin{equation}
\left(
\begin{array}{ccccc}
\frac12 & -1 & 1 & -1 &\ldots \\
\partial f_1/\partial t_0 & \partial f_1/\partial t_2 & \ldots \\
\partial f_2/\partial t_0 & \partial f_2/\partial t_2 & \ldots \\
\vdots \\
\partial f_{N-1}/\partial t_0 & \partial f_{N-1}/\partial t_2 & \ldots \\
\end{array}
\right)
\cdot
\left(
\begin{array}{ccc}
\Delta t_0\\
\Delta t_2\\
\Delta t_4\\
\vdots\\
\Delta t_{2N-2}\\
\end{array}
\right)
=
\left(
\begin{array}{ccc}
1-\sideset{}{'}\sum_n (-)^{n/2}t_n \\
-f_1\\
\vdots\\
-f_{N-1}\\
\end{array}
\right)
.
\label{eq.Ne}
\end{equation}
To keep track of the prime at the sum symbols, a binary symbol
which attains values of $2$ or $1$
is helpful:
\begin{defn} (Neumann symbol $\epsilon$)
\begin{equation}
\epsilon_0=2;\quad \epsilon_{>0}=1.
\end{equation}
\end{defn}
The two terms in (\ref{eq.fzer}) are denoted $f^{(a)}$ and $f^{(b)}$.
The second and higher rows in the matrix (\ref{eq.Ne}) are derivatives of
$f=f^{(a)}+f^{(b)}$, which are computed with the chain and multiplication rules.
The variable $\lambda$ is eliminated with the
aid of the derivative of (\ref{eq.g1}),
\begin{equation}
\frac{\partial \lambda }{\partial t_k} = -\frac{1}{\epsilon_k}
.
\end{equation}
\begin{rem}
This implementation is one variant out of many. The simplicity of the previous formula means
that the elimination of $\lambda$ and $\Delta \lambda$ from the pool of
unknowns produces no computational load.
\end{rem}
The derivatives of (\ref{eq.gx}) are
\begin{gather}
\frac{\partial g(\lambda x_j)}{\partial t_k} =
\frac{1}{\epsilon_k}\left(
T_k(\lambda^d x_j^d)-
x_j^d d\lambda^{d-1}\sum_{l\ge 1} t_lT_l'(\lambda^d x_j^d)
\right)
\label{eq.dglx}
,
\end{gather}
using the chain rule with the previous equation.
The first term in (\ref{eq.fzer}) is
\begin{multline}
f_j^{(a)}=\frac{t_0}{2}T_0\left[\left(\frac{t_0}{2}T_0(\lambda^d x_j^d)+t_1T_1(\lambda^d x_j^d)+\cdots\right)^d\right]
\\
+t_1T_1\left[\left(\frac{t_0}{2}T_0(\lambda^d x_j^d)+t_1T_1(\lambda^d x_j^d)+\cdots\right)^d\right]
\\
+t_2T_2\left[\left(\frac{t_0}{2}T_0(\lambda^d x_j^d)+t_1T_1(\lambda^d x_j^d)+\cdots\right)^d\right]
+\cdots
\label{eq.fa}
\end{multline}
with derivatives
\begin{gather}
\frac{\partial f_j^{(a)}}{\partial t_k}
=
\frac{1}{\epsilon_k}T_k[g^d(\lambda x_j)]
+
dg^{d-1}(\lambda x_j)\frac{\partial g(\lambda x_j)}{\partial t_k}
\sum_{l\ge 1} t_lT'_l[g^d(\lambda x_j)]
\label{eq.dfadt}
.
\end{gather}
The term at $l=0$ in the $l$-sum is
skipped since $T_0'(.)=0$.
(\ref{eq.dglx}) is inserted in front of the $l$-sum.

The second term in (\ref{eq.fzer}) is
\begin{equation}
f_j^{(b)}
=
\lambda \left(\frac{t_0}{2}T_0(x_j^d)+t_1T_1(x_j^d)+t_2T_2(x_j^d)+\cdots\right)
,
\label{eq.fb}
\end{equation}
with derivatives
\begin{gather}
\frac{\partial f_j^{(b)}}{\partial t_k} =
\frac{\partial \lambda }{\partial t_k}g(x_j)
+\lambda \frac{\partial g(x_j)}{\partial t_k}
=
-\frac{g(x_j)}{\epsilon_k}
+\lambda \frac{1}{\epsilon_k}T_k(x_j^d).
\label{eq.dfbdt}
\end{gather}
The sums of (\ref{eq.dfadt}) and (\ref{eq.dfbdt})
fill all but the first row of the matrix (\ref{eq.Ne}), and the negated sums
of (\ref{eq.fa}) and (\ref{eq.fb}) fill the right hand side.

\subsection{Convergence}
All digits of $t_n$ and $1/\lambda$ are considered stable which remain the same if
the order of the basis set $\{t_n\}$ is increased from $N$ to $N+4$, and these
stable digits will be shown in Section \ref{sec.res}\@.
The two (truncated) Chebyshev series are translated into the two
equivalent polynomials, and the stable (common) digits of $b_n$
are also reported. [This is an application of (\ref{eq.boft}) with
error propagation.] In this sense, the floating point representations of $t_n$
are $b_n$ are rounded towards zero.

\section{Results} \label{sec.res}
\subsection{z=2}
The most prominent solution is characterized by a leading order $z=2$ with
\begin{equation}
\lambda  \approx 0.39953.
\end{equation}
Broadhurst's 1018 most-significant digits of
$1/\lambda \approx 2.5029078$ are available via a reference
in the Online Encyclopedia of Integer Sequences \cite[A006891]{EIS},
updating an earlier value by Briggs \cite{BriggsMC57}.

The equivalent table of $n$ and $t_n$ at $d=1$ follows.
Long lines of digits are
wrapped around:

\scriptsize \begin{verbatim}
  0  0.5657908632724943155053234479753520221074493015133421406289563842809270230277394544246259
       011330003848845022918090496111898399
  2 -0.700391573973713787145980440332471261711550900968326558936525103514630069672843746266611
       2605202911962834490844304668109574045
  4  0.0173621867222441842092787099400084179514093184582945769724809006730767362263833677577478
       919456847233605679815143774768249724
  6  0.0006236559136940908585298153810290290296475164268085752829679471895290731944075421401943
       614274816952820691319610401477096491
  8 -0.000025266414376208310293595863274053777442956794506726672370626681912006533496105357064
       29129678242706445734230384800866127081
 10  0.0000002781260604292479399098215157081791712604067313923235379235837579995393225358243627
       3878484563185046160651019299165207065
 12  0.0000000077936819963404276117862621362965789982955679476321280713649199426441939964528134
       40475515594882948978268818766695877953
 14 -0.000000000327785722586730123873708759372184079817288268443710493022130004551717855508077
       34091127699425296427128013117312173163
 16  0.0000000000006384210078705431885183925429981891093422337470885619930147560493556137215402
       444436471651246493296658874924317760
 18  0.0000000000001778107114594587345498329317747650578350774121312001420048384738123759595693
       3133389845995766737758138453317984213
 20 -0.000000000000002799043155015243405083045051874952896840805205291021004938704453334036778
       277145181036915561331924809893972849572
 22 -0.000000000000000065987765706423742082420617124884422275834957432503771525513188833732520
       344678413898959592066716703455100407193
 24  0.0000000000000000031157462047608457677156560401194058182156452845428538122360790124887277
       888471555567698962440078037595332028476
 26 -0.000000000000000000026288817352565470028740677119624981235415929918927887971399803508489
       9487286649756918415632470426091310510068
 28 -0.000000000000000000001326155752923546645440506109645444391400468737885365499984259386868
       032271590021474125533154500539436762786
 30  0.0000000000000000000000456892132858964069886843715547017478914868661445083401985269517053
       73249571118425751022744438909374342921
 32 -0.000000000000000000000000261669143842729453902202384622846340198934739400123486545468197
       7688718781026782423507060753527153562631
 34 -0.000000000000000000000000019483941714088630622996200219342698291373135476338964862858461
       2208772721782152441705547209178529184368
 36  0.0000000000000000000000000005488618268001120285365797865347450259445755878549817589691912
       348942905406663752221458860366591795271
 38 -0.000000000000000000000000000000803615119146410423970695111862336617536776857636033883366
       5176549736209123317801265558837346131265
 40 -0.000000000000000000000000000000272928783224017159344918835131716069410115870190841420752
       8260064107651954953449841125396642705536
 42  0.0000000000000000000000000000000055802807643744090096015047177176805233416723117916922201
       603522647678299941997881614294928232295
 44  0.0000000000000000000000000000000000428579376944369989419311461634796486388452488691470855
       160050262628641836975040193582550337891
 46 -0.000000000000000000000000000000000003805534790674745785748165448883128965225994061896213
       8680306996844531622865057406490231737261
 48  0.0000000000000000000000000000000000000454730146146092569561537159214715109779392380848040
       024052991533878115732523670857341381746
 50  0.0000000000000000000000000000000000000012856880952631616613496264033316738376351888786159
       55772710048659056300864112109853599657
 52 -0.000000000000000000000000000000000000000048841173241110779531795133119845934307395101800
       77492098933787969566617466343321303425362
 54  0.0000000000000000000000000000000000000000001784169720275549045117234108914889596378406416
       6863250305019534622990073694338459577825
 56  0.0000000000000000000000000000000000000000000238344930590278575257139020978076800439813301
       6858895181885626215542830567042812186221
 58 -0.000000000000000000000000000000000000000000000533425905678224058971848422103488739756551
       31310220176451368406865863831223688863672
 60 -0.000000000000000000000000000000000000000000000003414909710470510411967146436853300862417
       751078813413975732929598114507412935572869
 62  0.0000000000000000000000000000000000000000000000003539028909832830183981073006466249273856
       4067409053977091623664707959698105127494
 64 -0.000000000000000000000000000000000000000000000000004485555444886975322920067691672933793
       194542250046135516632720727879043742021878
 66 -0.000000000000000000000000000000000000000000000000000113787377770035835911665832851820489
       23743586494443395196018665583594206702949
 68  0.0000000000000000000000000000000000000000000000000000045078100913470491750586844705707864
       42930571509190012568778013289412289813409
 70 -0.000000000000000000000000000000000000000000000000000000018133006944833701581226868346901
       64607157726409542422711423361367534302226
 72 -0.000000000000000000000000000000000000000000000000000000002174867119590406343851474044134
       81335237693943019828230475139169452124111
 74  0.0000000000000000000000000000000000000000000000000000000000496862104257630507989687593666
       374104723739661958257634148071260673007
 76  0.0000000000000000000000000000000000000000000000000000000000003050554484567852177812154431
       197959932218874902419033309048699799358
 78 -0.000000000000000000000000000000000000000000000000000000000000033498311561034693545843535
       0052156997519744655188887627916550362279
 80  0.0000000000000000000000000000000000000000000000000000000000000004527590094881243530685260
       5968118565891420964808271817046163737471
 82  0.0000000000000000000000000000000000000000000000000000000000000000103423598875438777825748
       5573036199162363130391497093308583693934
 84 -0.000000000000000000000000000000000000000000000000000000000000000000450299242051663454771
       795263762868632507049946111828727905335941
 86  0.0000000000000000000000000000000000000000000000000000000000000000000027428313773607458422
       7464939301109131771830302353395925019067
 88  0.0000000000000000000000000000000000000000000000000000000000000000000002009727142887944394
       2191940331186323737668132381183364540498
 90 -0.000000000000000000000000000000000000000000000000000000000000000000000005407644154626431
       9775699024893458068658018809505034372144
 92 -0.000000000000000000000000000000000000000000000000000000000000000000000000008467529148803
       04847834758872768166022849036943505097683
 94  0.0000000000000000000000000000000000000000000000000000000000000000000000000031985737186339
       3891418296413027668560883977465202924273
 96 -0.000000000000000000000000000000000000000000000000000000000000000000000000000057466843628
       4675810736663981678937691896929486580866
 98 -0.000000000000000000000000000000000000000000000000000000000000000000000000000000668154854
       46991658046723749549646398036058312243705
100  0.0000000000000000000000000000000000000000000000000000000000000000000000000000000451184387
       1270886465156982260768699951056609092057
102 -0.000000000000000000000000000000000000000000000000000000000000000000000000000000000510109
       3118547311281863095892199861517680972441
104 -0.000000000000000000000000000000000000000000000000000000000000000000000000000000000015095
       2899889893117535393879575008633927500623
106  0.0000000000000000000000000000000000000000000000000000000000000000000000000000000000005784
       8978366957498811062305983117450881881371
108 -0.000000000000000000000000000000000000000000000000000000000000000000000000000000000000002
       96860076774123816469785151254586608994118
110 -0.000000000000000000000000000000000000000000000000000000000000000000000000000000000000000
       26159630866409834309524824065939495168387
112  0.0000000000000000000000000000000000000000000000000000000000000000000000000000000000000000
       0675797965379667293610841149639913358407
114  0.0000000000000000000000000000000000000000000000000000000000000000000000000000000000000000
       0001185260268064102334926326440613443457
116 -0.000000000000000000000000000000000000000000000000000000000000000000000000000000000000000
       0000039471664898376094381214217934736500
118  0.0000000000000000000000000000000000000000000000000000000000000000000000000000000000000000
       0000007061272200790634490793605019948160
120  0.0000000000000000000000000000000000000000000000000000000000000000000000000000000000000000
       0000000078441087526952075341256121892401
122 -0.000000000000000000000000000000000000000000000000000000000000000000000000000000000000000
       00000000005391759440269970326885723499235
124  0.0000000000000000000000000000000000000000000000000000000000000000000000000000000000000000
       00000000000621433839232553661420859887097
126  0.0000000000000000000000000000000000000000000000000000000000000000000000000000000000000000
       0000000000001731657757671233420438125501
128 -0.000000000000000000000000000000000000000000000000000000000000000000000000000000000000000
       00000000000000067443060829753893881443456
130  0.0000000000000000000000000000000000000000000000000000000000000000000000000000000000000000
       00000000000000003652483640680563919080284
132  0.0000000000000000000000000000000000000000000000000000000000000000000000000000000000000000
       0000000000000000029557347361785303441624
134 -0.000000000000000000000000000000000000000000000000000000000000000000000000000000000000000
       00000000000000000000768611901159197870115
136 -0.000000000000000000000000000000000000000000000000000000000000000000000000000000000000000
       00000000000000000000001208200578812543434
138  0.0000000000000000000000000000000000000000000000000000000000000000000000000000000000000000
       00000000000000000000004393057915211811154
140 -0.000000000000000000000000000000000000000000000000000000000000000000000000000000000000000
       00000000000000000000000007775710345661442
142 -0.000000000000000000000000000000000000000000000000000000000000000000000000000000000000000
       00000000000000000000000000088568217982631
144  0.0000000000000000000000000000000000000000000000000000000000000000000000000000000000000000
       0000000000000000000000000005893642707841
146 -0.000000000000000000000000000000000000000000000000000000000000000000000000000000000000000
       0000000000000000000000000000006459531372031
148 -0.000000000000000000000000000000000000000000000000000000000000000000000000000000000000000
       000000000000000000000000000000019516562928
150  0.0000000000000000000000000000000000000000000000000000000000000000000000000000000000000000
       0000000000000000000000000000000071948895
152 -0.000000000000000000000000000000000000000000000000000000000000000000000000000000000000000
       00000000000000000000000000000000000311746
154 -0.000000000000000000000000000000000000000000000000000000000000000000000000000000000000000
       0000000000000000000000000000000000003291
156  0.0000000000000000000000000000000000000000000000000000000000000000000000000000000000000000
       00000000000000000000000000000000000007904
\end{verbatim}\normalsize

The equivalent stable estimates of the $b_n$ for this principal
solution at $z=2$
are:

\scriptsize \begin{verbatim}
  2 -1.527632997036301454035890310160257274974696028575842206144610860054595492964020577764313
       80701752738480877611758895682185543
  4  0.1048151947873037332167426138007218458656299992476409576574003078879413514272968522834696
       31229691364815580578035551868733
  6  0.0267056705251933540326520949436840400489828245712553869368102342195273415751956204909299
       55727292886661386456788292382
  8 -0.003527409660908709170234190769227694069753863863048995486587186788789137640735529902340
       66026364162201257631510784644
 10  0.0000816009665475317451721904864440992915367795680059157284744337637040509633793145598636
       7487595025088475730836602
 12  0.0000252850842339635361762625518763989120363749856139868813502169742225949626644722991135
       3537351897785177021222
 14 -0.000002556317166278493846353254082873397693909214123249051883015238706015469904131104513
       903982874702644269841
 16 -0.000000096512715508912032163725767577852903309294329912746415412354641940246862574197126
       7191899550777987936
 18  0.0000000281934639745040913707566272517012381087057565077399072734902349656431332938635785
       65042561426215514
 20 -0.000000000277305116079901172437311261354608894542298492401176729164796655006915246627235
       5751955678092972
 22 -0.000000000302842702213056632983880872952233308909357930141156404683591761118364361653642
       23366616410882
 24  0.0000000000267058928074807555396472475218453020927141066761240072926251938510138459085654
       76029894011
 26  0.0000000000009962291641028482309516423739168010973641395908612318394179172812550254542680
       6722454314
 28 -0.000000000000362420298290415608223127234689812443132148566311172785414580761787322163007
       369943524
 30  0.0000000000000217965774482707043220507000979414607626330551237338227994433136351844298928
       89151245
 32  0.0000000000000015292328994809633812036270843646278271962252136334146277442594850713415288
       753536
 34 -0.000000000000000318472878995278789748224049950902399205080342381778791026825878299688402
       28761
 36  0.0000000000000000113467210621204186902773212125499882109020286136044079317046380395423364
       80
 38  0.0000000000000000018816760568232668841203920682045414393478508433015342382297496601166865
       1
 40 -0.000000000000000000227561256460702582887704617785944487885375708180060065740481154582306
 42 -0.00000000000000000000098224476556928294515508544402787045557333578831363282612181126123
 44  0.0000000000000000000020641297584976918664437582373005112128267581989566512117513015412
 46 -0.000000000000000000000124932008068648399476793391930454469740646818703722594463270599
 48 -0.0000000000000000000000107706103770474744566181106412701549392134721586703334688521
 50  0.0000000000000000000000018727456323338913159319680190999726028324097866404744001097
 52 -0.000000000000000000000000025776334884038022483024846450586376516821232359223724208
 54 -0.00000000000000000000000001554233189736463104356497580604304978763137037382869610
 56  0.0000000000000000000000000012806626551535968818106767670009804575977590900839329
 58  0.000000000000000000000000000055750205796890611189266461733394520137627182303249
 60 -0.00000000000000000000000000001523753513286678958880994208647025914164300165148
 62  0.000000000000000000000000000000489553018633832353824639228798471893324776739
 64  0.000000000000000000000000000000106224525652671030953252135604394033249927799
 66 -0.00000000000000000000000000000001129930672170938039709583095905246032244541
 68 -0.00000000000000000000000000000000024459998919548209728141752362355767054563
 70  0.000000000000000000000000000000000120363378062993272576012717610146326113
 72 -0.00000000000000000000000000000000000563012978201117079658608329257623690
 74 -0.00000000000000000000000000000000000075237405185431847091563898005564015
 76  0.0000000000000000000000000000000000001034729036462367522938060502393519
 78  0.0000000000000000000000000000000000000000559027103251345569985604893845
 80 -0.000000000000000000000000000000000000000980444278474974044174175006208
 82  0.000000000000000000000000000000000000000067458101888888796308753858680
 84  0.00000000000000000000000000000000000000000454649630780742528187568226
 86 -0.00000000000000000000000000000000000000000095809721960822642394780771
 88  0.0000000000000000000000000000000000000000000251772819390575385916219
 90  0.000000000000000000000000000000000000000000007227639070634795496186
 92 -0.000000000000000000000000000000000000000000000756392738994815737804
 94 -0.00000000000000000000000000000000000000000000001316892955819095928
 96  0.00000000000000000000000000000000000000000000000802592842968755113
 98 -0.00000000000000000000000000000000000000000000000045956332222745714
100 -0.00000000000000000000000000000000000000000000000004121498443984379
102  0.0000000000000000000000000000000000000000000000000072557224505846
104 -0.0000000000000000000000000000000000000000000000000001465310222330
106 -0.000000000000000000000000000000000000000000000000000056265758540
108  0.00000000000000000000000000000000000000000000000000000543920805
110  0.000000000000000000000000000000000000000000000000000000121289603
112 -0.000000000000000000000000000000000000000000000000000000058647890
114  0.000000000000000000000000000000000000000000000000000000003149823
116  0.000000000000000000000000000000000000000000000000000000000306407
118 -0.000000000000000000000000000000000000000000000000000000000050866
120  0.000000000000000000000000000000000000000000000000000000000000898
\end{verbatim}\normalsize

This corrects up to five digits in some rows of
Lanford's table, where error bars are of the order $10^{-30}$ \cite{LanfordBAMS6}.

\subsection{Extra z=2}
\begin{figure}
\includegraphics[scale=0.6]{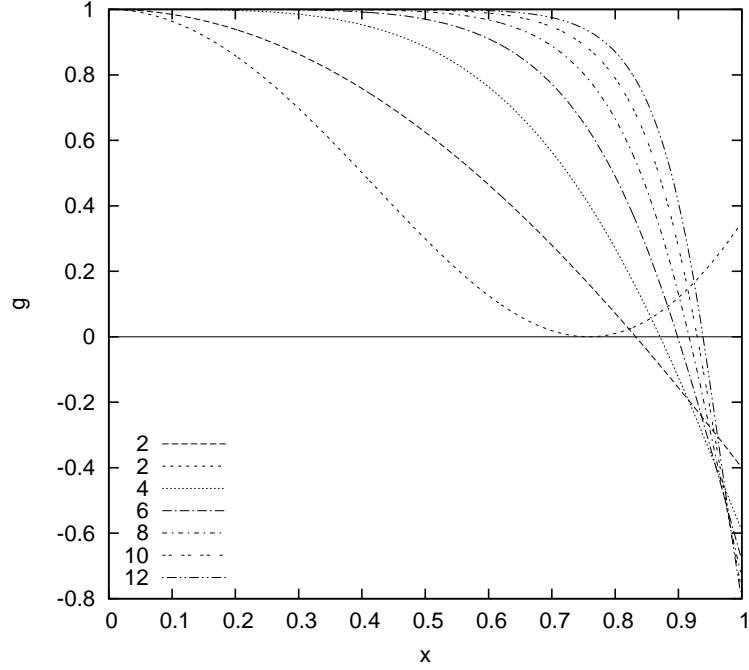}
\caption{
Solutions $g(x)$ discussed in this work, marked with their $z$-parameters.
The three curves at $z=2$ or $z=4$ have been plotted
by Campanino et.\ al.\ \cite{CampaninoTop21}.
}
\label{fig.1}
\end{figure}

The solution recorded by Stephenson and Wang at $z=2$ which returns
to a positive value $g(1)>0$ (see Figure \ref{fig.1}), hence assumes a negative value of
$\lambda$, is also refined
\cite{StephensonAML4,StephensonAML4b}.

Its table of $n$ and $t_n$, value of $1/\lambda$, and table
of $n$ and $b_n$ are:
\scriptsize \begin{verbatim}
  0  0.695239872293346037618864099266551483193696717452792702
  2 -0.2941925542257057866921887876571014540207923552631288627
  4  0.33028072470630568964000834403239424732161343709604632350
  6 -0.0315368919233901472203850867002152965050184294468376104
  8 -0.0028714832609075151818101565163240137621094990590925332
 10  0.00073699063183551803763665511752866415092050599857134560
 12 -0.0000289143566935215723183757294107103328257404482568374
 14 -0.00000639269493477664368809268119272428349626603668195306
 16  0.00000089250734493829604749070500496361151531437475099399
 18 -0.000000005841530325386329625967647262112183068522122098459
 20 -0.00000000904338569609666460923438628990512000571101803815
 22  0.000000000794308476080916066911148179642789312890062120789
 24  0.00000000003025917767192722188179590965141382923724510004
 26 -0.0000000000105467686794658509993062892575928743129899731754
 28  0.0000000000005178028438818918034562105519025020884646961949
 30  0.000000000000067093995175319688853016991708074433815242496
 32 -0.000000000000010559545218086633148414702577521791670273535
 34  0.000000000000000154347918061355616738318944738955043641955
 36  0.000000000000000094274699599167949022401989148613318504987
 38 -0.000000000000000009033368564852075748116667512849206793724
 40 -0.00000000000000000022125271892120189926978367638849260329
 42  0.000000000000000000106886427363051040686657229276492155630
 44 -0.0000000000000000000062817362134522292817164837731424606819
 46 -0.0000000000000000000005493232753850121642056491695464972991
 48  0.0000000000000000000001043062141037156724808203246971504313
 50 -0.00000000000000000000000285035665884706861429840914755099124
 52 -0.00000000000000000000000078515923793851604089335622069042986
 54  0.00000000000000000000000008869529303851488865811270044693548
 56  0.000000000000000000000000000649860120188466719127047774058709
 58 -0.00000000000000000000000000090206496864462618542641659904916
 60  0.00000000000000000000000000006377438397374861263138389614442
 62  0.0000000000000000000000000000036876697756045905245014394098439
 64 -0.00000000000000000000000000000089470859900052789175921501907
 66  0.00000000000000000000000000000003402948808776406735817735563
 68  0.00000000000000000000000000000000589033257732122101931819868
 70 -0.0000000000000000000000000000000007780486465952083239308934
 72  0.00000000000000000000000000000000000403483098534352086273508
 74  0.00000000000000000000000000000000000706836056123442192425924
 76 -0.00000000000000000000000000000000000058180447831663676039050
 78 -0.0000000000000000000000000000000000000222072831728450630950
 80  0.0000000000000000000000000000000000000072079140255579564551
 82 -0.0000000000000000000000000000000000000003430548564035407762
 84 -0.00000000000000000000000000000000000000004178999679170301237
 86  0.0000000000000000000000000000000000000000064420879852046845
 88 -0.0000000000000000000000000000000000000000000991742129585092
 90 -0.00000000000000000000000000000000000000000005317441241661680
 92  0.0000000000000000000000000000000000000000000050075104582881
 94  0.0000000000000000000000000000000000000000000001177289624901
 96 -0.0000000000000000000000000000000000000000000000561911232050
 98  0.0000000000000000000000000000000000000000000000031935471013
100  0.0000000000000000000000000000000000000000000000002844395171
102 -0.0000000000000000000000000000000000000000000000000518262434
104  0.0000000000000000000000000000000000000000000000000012926078
106  0.0000000000000000000000000000000000000000000000000003882830
108 -0.0000000000000000000000000000000000000000000000000000418767
110 -0.00000000000000000000000000000000000000000000000000000044
\end{verbatim}\normalsize

\scriptsize \begin{verbatim}
F = -2.857124135141400000343125136089134962070298108661379658
  2 -3.6682158140188213295544406557233531750716131334358113
  4  3.3898762669371844947446408632890174631221528494613
  6  0.57625523018561194477776150019529828845365167129
  8 -1.267984110033081623678420832380008814992950983
 10  0.1604899515138390776344657544573363292210139
 12  0.26867432609164832115329089526627621696060
 14 -0.100525526259689645031469710302302791466
 16 -0.0318348172556850038010772555453835216
 18  0.026218359709652495381298970181341915
 20  0.00087630890012302809197295880627673
 22 -0.005064121738585709401711490844130
 24  0.000726562430533655974197555993
 26  0.0008382684112571983909288292855
 28 -0.000291516068557741383544793566
 30 -0.0001071802944726395309771115
 32  0.00007685496751041026904906
 34  0.00000504836963733199948456
 36 -0.0000154135179534592504525
 38  0.000002187336070281942600
 40  0.00000236018113386321604
 42 -0.0000008586990279492577
 44 -0.000000253360946249660
 46  0.000000200634693258443
 48  0.00000000691689450040
 50 -0.00000003682574116969
 52  0.0000000056916964608
 54  0.0000000054513684210
 56 -0.000000002036725611
 58 -0.00000000057562842
 60  0.00000000046692491
 62  0.00000000001246505
 64 -0.0000000000837368
 66  0.0000000000134892
 68  0.000000000012016
 70 -0.0000000000045542
 72 -0.000000000001245
 74  0.000000000001017
\end{verbatim}\normalsize

\subsection{z=4}
A known solution with smoothness $z=4$ is accurately described by
the following table of $n$ and $t_n$ at $d=2$
\cite{StephensonAML4,StephensonAML4b}:

\scriptsize \begin{verbatim}
  0  0.3259810020643102220511536428215443201428992452133
  2 -0.800208447142190906917818576407006363536002907994
  4  0.0418898841904523139724543647479542658619459241985
  6  0.0043894976258619771863959729197546595065618948810
  8 -0.00069010337244814876166850830661956401216977657241
 10  0.0000144814075746787149177686108721351459836488957
 12  0.00000476300216712418243311831014922536070727006694
 14 -0.000000494555259455886663361734329319366255446136691
 16 -0.000000002894498873151192531198855036171513344547245
 18  0.0000000043659575703295075770575545118494173777559517
 20 -0.000000000278893954012460398480375551751458327306229
 22 -0.00000000001905793842739537879484162909974046955218144
 24  0.00000000000387515529016371762223417515171078981195
 26 -0.00000000000010811170201073782340784343123593525691
 28 -0.00000000000002781706617592988574104718860817240835
 30  0.000000000000003140943030258792107754156277223196029
 32  0.0000000000000000113037783599017151502130939738831523
 34 -0.0000000000000000294068162559927660239320286755977003
 36  0.00000000000000000217380400301654878938939813471539029
 38  0.00000000000000000009589354775830741278121051671184423
 40 -0.000000000000000000026674587926268032143341253347671646
 42  0.0000000000000000000011684247116706485719678443835458607
 44  0.00000000000000000000014912430927612388099838950079617660
 46 -0.00000000000000000000002158773018436811136342959944383313
 48  0.00000000000000000000000028451466280962684221628155260692
 50  0.000000000000000000000000171842585248494306526822030351909
 52 -0.00000000000000000000000001539313006678710338179573132491
 54 -0.00000000000000000000000000039226023013920647068619398834
 56  0.0000000000000000000000000001677001372390369884707804569
 58 -0.0000000000000000000000000000090135091823733200907687260
 60 -0.0000000000000000000000000000008349267645798270024985074
 62  0.0000000000000000000000000000001434920035967166259561752
 64 -0.00000000000000000000000000000000319516654856190813929
 66 -0.00000000000000000000000000000000105418754132973448968
 68  0.000000000000000000000000000000000107373483256963049299
 70  0.00000000000000000000000000000000000151525141050310900
 72 -0.0000000000000000000000000000000000010833036852367435
 74  0.000000000000000000000000000000000000067025066870316
 76  0.0000000000000000000000000000000000000048204263164179
 78 -0.0000000000000000000000000000000000000009684033207264
 80  0.000000000000000000000000000000000000000028596298782
 82  0.0000000000000000000000000000000000000000066708949798
 84 -0.00000000000000000000000000000000000000000076077044388
 86 -0.0000000000000000000000000000000000000000000036628899
 88  0.00000000000000000000000000000000000000000000721368704
 90 -0.0000000000000000000000000000000000000000000005098973
 92 -0.0000000000000000000000000000000000000000000000273795
 94  0.00000000000000000000000000000000000000000000000672075
 96 -0.000000000000000000000000000000000000000000000000257959
 98 -0.0000000000000000000000000000000000000000000000000418300
100  0.000000000000000000000000000000000000000000000000005520
\end{verbatim}\normalsize

The associated Feigenbaum constant $1/\lambda$ (updating the value
of $1.690302$ \cite{BriggsJPA24}) and equivalent
table of $n$ and Taylor coefficients
$b_n$ are:

\scriptsize \begin{verbatim}
F = 1.690302971405244853343780150324161348228278059709
  4 -1.83410790700941066477722032786167658753580656671
  8  0.012962226191371748194249954526500692423364159
 12  0.3119017366428453740938210685407183598371022
 16 -0.0620146232838494154168020915170780000963
 20 -0.037539476018044801855283971780562014903
 24  0.0176482141699045721668837916436744792
 28  0.00193502990250861524539045696343944
 32 -0.002811394115124136184245956801096
 36  0.00009519227150370368772733484141
 40  0.000435491310822346046233354127
 44 -0.000075173146500397955569561474
 48 -0.00006736728882241037232621045
 52  0.0000269344724659457260827140
 56  0.000006286772295217247380353
 60 -0.0000059502454648295960937
 64  0.0000002101623794439990936
 68  0.0000009189230440398699384
 72 -0.00000021302788365308075
 76 -0.00000010529385239013139
 80  0.0000000526461946904
 84  0.0000000076135782693113
 88 -0.000000009844458690011
 92  0.00000000047080672420
 96  0.00000000154345662656
100 -0.00000000035684794770
104 -0.00000000018996907837
108  0.00000000009357108066
112  0.000000000013961433
116 -0.000000000017731035
120  0.00000000000087943
124  0.00000000000272384
128 -0.0000000000006104
\end{verbatim}\normalsize

\subsection{z=6}

The Chebyshev coefficients $t_n$ for $z=6$ with $d=3$ are:

\scriptsize \begin{verbatim}
  0  0.21021341309304017470311495545314716311007687224938007663559
  2 -0.85061644791828709112050514504136245282682454755905389538579
  4  0.05628161905542202323276948575124712239104414802830965331397
  6  0.009935058137903504135993082540678516096531085856931659844078
  8 -0.00208555090028753720240399596899815942535546644082405772236
 10  0.00002667282124614707868884031418105463714552279637627290637
 12  0.00003829372515583810559086517586800603514191708130147980083
 14 -0.000004625699664422667849249798381290140479802593287643188046
 16 -0.0000002831151594281039534174488537456346214130161032444781955
 18  0.00000012338683107538209784038150537174370653614397236744514266
 20 -0.0000000068174045281617081780530926809083421690600186685820884
 22 -0.000000001899873267576372584804450734908402877469178801808811196
 24  0.000000000346392675995598419474794310398532136946869032638576
 26  0.00000000000390068105896541570853177397127706284297914293009653
 28 -0.0000000000079149906752964255297343374974010828294089885968213
 30  0.00000000000077931694780771578707264938426098858651135113840855
 32  0.0000000000000852773790825838513873257851533380830018759450045
 34 -0.0000000000000260383569822231813226762772106465631988119103374
 36  0.0000000000000010025050056296133091452636511043678546651959869
 38  0.000000000000000449476081511898302858473633434434053219462678
 40 -0.0000000000000000700314136799831204483340156830840996622492831
 42 -0.00000000000000000211427500105270276952688421271708844343409553
 44  0.0000000000000000017355760875149987461939747544759876131044824
 46 -0.000000000000000000141052513634382469967859280854950424746015099
 48 -0.000000000000000000022360690323524459052845995488637640910740059
 50  0.000000000000000000005463391250275824699752989313602943626365874
 52 -0.0000000000000000000001032678941939516175678341770001697926196876
 54 -0.00000000000000000000010759324343630651925920681050043268884835996
 56  0.000000000000000000000013852158702642782628823590465480410418581454
 58  0.0000000000000000000000008484209644875932810101724300513674300141197
 60 -0.00000000000000000000000039364298268040750575103277016650317751108005
 62  0.000000000000000000000000024020661589849142360039060864756583132179518
 64  0.000000000000000000000000006001745276715558023520845370409650872133890
 66 -0.0000000000000000000000000011746748308382063766464342078813559899538285
 68 -0.0000000000000000000000000000041146309997350933302276992059245312709099
 70  0.0000000000000000000000000000261692711387768568430763903240667516315367
 72 -0.0000000000000000000000000000027486585506723445468000794116756014447496
 74 -0.000000000000000000000000000000275308623842612943860805940900180953217
 76  0.000000000000000000000000000000090196778810492963298456791888344590739
 78 -0.00000000000000000000000000000000376607611498847643758409590315035925
 80 -0.0000000000000000000000000000000015696826935725241130869385256205437
 82  0.00000000000000000000000000000000025350782115307816392785281062157
 84  0.000000000000000000000000000000000006853080930004380405185843345
 86 -0.00000000000000000000000000000000000631173846425375219150669720605
 88  0.0000000000000000000000000000000000005372233333770510172444422828
 90  0.00000000000000000000000000000000000008002743541502796444109333872
 92 -0.00000000000000000000000000000000000002051613050532287398283675699
 94  0.0000000000000000000000000000000000000004791508300478340494151358
 96  0.000000000000000000000000000000000000000398185100966144593387684
 98 -0.000000000000000000000000000000000000000053917069856953084333496
100 -0.0000000000000000000000000000000000000000028656289415773532604205
102  0.000000000000000000000000000000000000000001495957563527673051760
104 -0.000000000000000000000000000000000000000000100085074963152324806
106 -0.00000000000000000000000000000000000000000002193992906351529577
108  0.00000000000000000000000000000000000000000000458733639392039182
110 -0.000000000000000000000000000000000000000000000017873394641713124
112 -0.000000000000000000000000000000000000000000000098786806004482188
114  0.00000000000000000000000000000000000000000000001114827380233477
116  0.000000000000000000000000000000000000000000000000956056408698808
118 -0.000000000000000000000000000000000000000000000000348947364249635
120  0.000000000000000000000000000000000000000000000000016973874478220
122  0.00000000000000000000000000000000000000000000000000581996626627
124 -0.000000000000000000000000000000000000000000000000001006284719313
126 -0.00000000000000000000000000000000000000000000000000001843151562
128  0.000000000000000000000000000000000000000000000000000024126102022
130 -0.000000000000000000000000000000000000000000000000000002219387394
132 -0.000000000000000000000000000000000000000000000000000000289110236
134  0.000000000000000000000000000000000000000000000000000000080263669
136 -0.0000000000000000000000000000000000000000000000000000000023710691
138 -0.000000000000000000000000000000000000000000000000000000001508066
140  0.0000000000000000000000000000000000000000000000000000000002162263
142  0.0000000000000000000000000000000000000000000000000000000000096712
144 -0.00000000000000000000000000000000000000000000000000000000000581212
146  0.0000000000000000000000000000000000000000000000000000000000004193
148  0.00000000000000000000000000000000000000000000000000000000000008236
150 -0.000000000000000000000000000000000000000000000000000000000000018198
152  0.000000000000000000000000000000000000000000000000000000000000000165
\end{verbatim}\normalsize

The Feigenbaum constant (updating $1.4677$ \cite{BriggsJPA24})
and the Taylor coefficients $b_n$ are:
\scriptsize \begin{verbatim}
F = 1.4677424503199009444538343151089737463687971293967614336021
  6 -1.907736962634173070120370338436410222425918058369139834419
 12 -0.3328834745034284193531634790623119500310257123750925067
 18  0.7127013698436617555865888754368253815367371119995311
 24  0.03518312285829711535877045080703530781998304860291
 30 -0.27250710400431126009882590314277478262909841416
 36  0.025924899594891682578009471584373468628157938
 42  0.0935345652691494827767146604329617762239436
 48 -0.015200629738589197533218095454339549332587
 54 -0.0361996943354910461489660634048101524283
 60  0.00842218588897031494780746155241064521
 66  0.0152369615087910404306655219832892611
 72 -0.00568137035329806459395752180194685
 78 -0.0056690875578583361086105903358374
 84  0.00320118052966879223649858259047
 90  0.0018209196254639367945467799981
 96 -0.001494614499707863586495267381
102 -0.00055268652092897063886704420
108  0.000651366223545191258668000
114  0.00016211183897711159615729
120 -0.0002832959350896794399643
126 -0.000039058079635634421561
132  0.00012203126950356924900
138  0.0000032127601587875791
144 -0.000050975296723247046
150  0.000004060585842767806
156  0.00002068636917088808
162 -0.0000037601470720777
168 -0.000008262127271535
174  0.00000239087934851
180  0.0000032522821183
186 -0.000001347046590
192 -0.000001242833132
198  0.000000708550725
204  0.00000045111756
210 -0.0000003523916
216 -0.000000151591
222  0.000000166981
228  0.000000044980
234 -0.00000007591
240 -0.0000000101
\end{verbatim}\normalsize

\subsection{z=8}
The Chebyshev coefficients $t_n$ for $z=8$ with $d=4$ are:

\scriptsize \begin{verbatim}
  0  0.1381286283736299578601368862
  2 -0.8842303740823175512755608808
  4  0.0665209634842108260820275058
  6  0.0160708249063224794092104569
  8 -0.0038907994146791521600746136
 10 -0.0000097551920092753105368942
 12  0.000123849745247201425628817
 14 -0.00001527270699745541145444107
 16 -0.0000021352008988295402412375489
 18  0.000000772398623891616273984173
 20 -0.000000024219287650319566438404
 22 -0.000000023416621864635256604937
 24  0.000000003826707490135532768257
 26  0.00000000032042363484105018598
 28 -0.00000000017467929304630940850
 30  0.00000000001152216008656868781
 32  0.000000000004587462209785440857
 34 -0.000000000000995958581659698440
 36 -0.0000000000000312775335815343648
 38  0.000000000000038873506228181227
 40 -0.00000000000000392306669959245007
 42 -0.00000000000000085951077958354364
 44  0.000000000000000249251998355112348
 46 -0.0000000000000000018980520249321427
 48 -0.00000000000000000856801506577853604
 50  0.000000000000000001177536979181658135
 52  0.0000000000000000001511172911297034026
 54 -0.00000000000000000006142892999320489
 56  0.00000000000000000000261014286449232
 58  0.0000000000000000000018461512874019
 60 -0.000000000000000000000333176674195
 62 -0.000000000000000000000022406130989
 64  0.000000000000000000000014800851697
 66 -0.00000000000000000000000113012074
 68 -0.00000000000000000000000038099450
 70  0.00000000000000000000000009008411
 72  0.0000000000000000000000000018534
 74 -0.00000000000000000000000000346207
 76  0.0000000000000000000000000003859
 78  0.000000000000000000000000000073
\end{verbatim}\normalsize

The corresponding Feigenbaum constant
(updating $1.35798$ \cite{BriggsJPA24}) and the associated
Taylor coefficients $b_n$ are:
\scriptsize \begin{verbatim}
F = 1.3580172791380503454873763331
  8 -1.89735300467491340532977247
 16 -0.73884380388552253531567
 24  0.989774470130239041388
 32  0.445857788637156558
 40 -0.5879109210650583
 48 -0.268029707317435
 56  0.326144611238
 64  0.20652701084
 72 -0.201257882
 80 -0.1608771
 88  0.1402826
 96  0.1105
104 -0.094
\end{verbatim}\normalsize

\subsection{z=10}
The Chebyshev coefficients $t_n$ for $z=10$ with $d=5$ are:

\scriptsize \begin{verbatim}
  0  0.087844965370486707
  2 -0.9093012275250430520
  4  0.074588019670768413
  6  0.022306567620926918
  8 -0.00591629904025435163
 10 -0.00011666771033891637
 12  0.000271811762145798921
 14 -0.0000328109239789392
 16 -0.00000770036081694509
 18  0.00000258451857056376
 20 -0.000000003534676935799
 22 -0.000000121746886087767
 24  0.0000000176227062237214
 26  0.00000000328463293938500
 28 -0.00000000132079164712664
 30  0.000000000033312229401120
 32  0.000000000058825803651164
 34 -0.00000000001012188314110
 36 -0.0000000000013849849099
 38  0.00000000000069207426816
 40 -0.0000000000000327774239
 42 -0.0000000000000290338777
 44  0.0000000000000058170125
 46  0.000000000000000577884
 48 -0.000000000000000369611
 50  0.000000000000000025399
 52  0.00000000000000001449192
 54 -0.0000000000000000033648
 56 -0.000000000000000000223
\end{verbatim}\normalsize

This translates into $1/\lambda$ and the Taylor coefficients $b_n$
as follows:
\scriptsize \begin{verbatim}
F = 1.2915168672623445696
 10 -1.8517140134795485
 20 -1.124743004799023
 30  1.0746332407420
 40  1.07554234877
 50 -0.6752703771
 60 -0.99566771
 70  0.2886430
 80  0.99702
 90 -0.0150
\end{verbatim}\normalsize

\subsection{z=12}
The list of $n$ and $t_n$ for $z=12$ with $d=6$ starts:

\scriptsize \begin{verbatim}
  0  0.050315274170474025
  2 -0.9292027766462839346
  4  0.081321549028015737
  6  0.0284367418795728353
  8 -0.00805828726265292218
 10 -0.00029782425094054862
 12  0.000483912632074647391
 14 -0.00005617186728366212
 16 -0.00001926578700820676
 18  0.000006201794732562265
 20  0.000000194426753767854
 22 -0.0000004017148408827121
 24  0.00000005067173288804275
 26  0.00000001645633128370378
 28 -0.000000005657763568442150
 30 -0.0000000001193549987650470
 32  0.000000000365299595939128
 34 -0.000000000049892905014794
 36 -0.000000000014698792552832
 38  0.00000000000539147386315
 40  0.00000000000006412134576
 42 -0.00000000000034507823335
 44  0.00000000000005035434125
 46  0.0000000000000135831268
 48 -0.0000000000000052890573
 50 -0.0000000000000000142316
 52  0.00000000000000033350259
 54 -0.0000000000000000518608
 56 -0.00000000000000001272412
 58  0.0000000000000000052730
 60 -0.00000000000000000003532
 62 -0.0000000000000000003262
 64  0.0000000000000000000539
\end{verbatim}\normalsize

This translates into $1/\lambda$ and Taylor coefficients $b_n$:
\scriptsize \begin{verbatim}
F = 1.2465277517207492954
 12 -1.79116162311222203
 24 -1.463168537263831
 36  0.9856559520745
 48  1.76384751144
 60 -0.352844894
 72 -1.91762338
 84 -0.54023
 96  2.0135
108  1.521
120 -1.866
132 -2.43
144  1.
\end{verbatim}\normalsize

\subsection{z=14}
The list of $n$ and $t_n$ for $z=14$ with $d=7$ starts:

\scriptsize \begin{verbatim}
  0  0.0210003942667285
  2 -0.945643213077312
  4  0.08714605274152282
  6  0.03437323735545195
  8 -0.01025712295409548
 10 -0.000548744504139
 12  0.000757490037990091
 14 -0.000083828165279389
 16 -0.000038842635955055
 18  0.000012160982312532
 20  0.0000007928961871139
 22 -0.000001006742945403
 24  0.00000010854938696682
 26  0.00000005550666127079
 28 -0.000000017168885225179
 30 -0.0000000012329704551329
 32  0.0000000014734835693269
 34 -0.0000000001538909310688
 36 -0.000000000083949894876
 38  0.000000000025491091508
 40  0.000000000002018718384
 42 -0.00000000000224862047
 44  0.00000000000022506633
 46  0.0000000000001314485
 48 -0.00000000000003902405
 50 -0.00000000000000336485
 52  0.0000000000000035144
 54 -0.00000000000000033657
 56 -0.000000000000000209532
 58  0.00000000000000006079
 60  0.00000000000000000565
\end{verbatim}\normalsize

This translates into the Feigenbaum constant and Taylor coefficients $b_n$:
\scriptsize \begin{verbatim}
F = 1.21391238764424391
 14 -1.72516768360581
 28 -1.7485120998868
 42  0.76642269718
 56  2.385250550
 70  0.383241
 84 -2.65092
 98 -2.3006
112  2.30
126  4.6
\end{verbatim}\normalsize

\section{Summary}
Three representations of the Feigenbaum Function $g(x^z)$
for orders
$z=2$ and $z=4$ have been computed with higher
precision than previously published.
The principal solutions in the parameter range $z=6$--$14$, some of which
have been characterized
in the literature by Feigenbaum constants with 5-digit accuracy, have been
made explicit.

\bibliographystyle{amsplain}
\bibliography{all}

\end{document}